
\magnification=1200  \hoffset.8truein \voffset 1truein 
\hsize30truepc \vsize50truepc  \let\bls\baselineskip  \parskip=.3\bls

\font\text=cmr10 \text \font\title=cmr12 \let\bb\relax 
\font\abst=cmr9  \font\SMC=cmr8  \def\smc #1{\SMC{#1}}    

\def\St{\S\thinspace}  \let\tsp\thinspace
\def\dAsh#1{\nobreak\tsp#1\penalty0\tsp}  \def\Space{\hbox{\space}}
\def\dash{\dAsh{--}} \def\Dash{\dAsh{---}}  \def\ie{i.\tsp e.} 
\def\words#1{\hbox{\quad #1\quad}} \def\Quad{\hbox{\quad}} 
\let\vep\varepsilon \let\vph\varphi \def\div{\,|\,}  \let\obar\overline

\def\CTRL#1{\centerline{\title #1}}  \def\Ctrl#1{\centerline{\smc #1}}

\CTRL{Factorization of integers and arithmetic functions} \vskip 1\bls
\Ctrl{LINCOLN DURST} \vskip2\bls 

\bgroup\narrower\narrower\abst\parindent=0pt \baselineskip=11truept
\centerline{ABSTRACT}

Elementary proofs of unique factorization in rings of arithmetic
functions using a simple variant of Euclid's proof for the fundamental
theorem of arithmetic.

\egroup\vskip.5\bls 

\vskip.5\bls{\noindent\bf1. Introduction.}

In {\sl The Elements} [{\bf11}, {\smc BOOKS VII} and {\smc IX}] Euclid
proved that, except for the order in which the factors are written,
positive integers greater than $1$ can be expressed uniquely as a
product of {\it irreducible\/} integers $p$.

Euclid's argument has two main components:
\item{(1)} Each {\it composite\/} integer $ab$, with $a>1$, $b>1$, 
has an irreducible factor $p$, \ie, a factor whose only divisors are
$1$ and $p$.  [{\smc BOOK VII, PROPOSITIONS \smc 31, 32}]
\item{(2)} Every irreducible integer is a {\it prime}, \ie, if $p$ 
is irreducible, $ab>1$, and $p\div ab$, then $p\div a$ or $p\div
b$. [{\smc VII, 20} and {\smc IX, 14}]

Euclid or Euclid's translators introduced the expression prime number
to represent what may be referred to in terminology of twentieth
century algebra as an irreducible positive integer.  Euclid deserves
full credit, of course, for recognizing that irreducible positive
integers have the more important property described in (2).

The need for a distinction between irreducibles and primes arose early
in the nineteenth century when an impasse was encountered while
attempting to prove Fermat's Last Theorem, since Euclid's (2) fails to
hold in some rings where (1) holds, especially in cases involving real
or complex algebraic numbers.  See Edwards [{\bf10}], Hardy and Wright
[{\bf12}], and Jacobson [{\bf14}].  

\noindent{\bf2. Euclid's original arguments.} 

For Euclid all integers are positive\Dash and ordered, indeed well
ordered, meaning no sequence of positive integers in which each member
is less than its predecessor can be an infinite sequence, but must
terminate after finitely many terms.

Here we are interested in factorization in a family of rings due to
Liouville that are not ordered.  First we modify Euclid's method,
using simple algebraic concepts, to reduce its dependence on the order
relation Euclid used.  Our variant of Euclid's (1) allows us to
preserve Euclid's (2).  Unlike the Fermat case, our procedure allows
us to use only rational integers.  In this sense our argument remains
entirely elementary.

Euclid begins {\smc BOOK VII} with a list of definitions, the relevant
part here being that numbers larger than one must be {\it
irreducible\/} or {\it composite}.  The former have only trivial,
{\ie, unit}, proper divisors; the latter have nonunit proper divisors.

Euclid's argument for property (1): If $a$ and $b$ are more than one,
$ab$ is composite.  If $a$ is irreducible, then $ab$ has an
irreducible factor.  Otherwise $a$ is composite, so some $c$ divides
$a$.  If $c$ is composite we repeat the previous step and begin
generating a chain of proper divisors of $ab$.  We must encounter an
irreducible factor, for the only alternative is to generate an
infinite sequence of divisors, each a proper divisor of its
predecessor, and that is impossible.  [{\bf11}, vol. 2, p. 332]

Euclid has two proofs for (2): {\smc FIRST.} If $p\div ab$, then
$pc=ab$ and $p:b=a:c$, \ie, $p/b=a/c$. If $p$ does not divide $b$, the
only common divisor of $p$ and $b$ is $1$, so the ratio $p:b$ is in
lowest terms and $m\ge1$ exists with $a=mp$, $c=mb$, hence $p\div
a$. [{\bf11}, vol. 2, p. 321]

{\smc SECOND.} Let $a$ be the least number divisible by all the irreducibles
$p_1, p_2, \dots, p_k$.  Suppose $q$ is an irreducible distinct from
each of $p_1, p_2, \dots, p_k$, and assume $q$ divides $a$.  Then
$a=qm$.  Since $q$ does not equal any of $p_1, p_2\dots, p_k$,
all the $p_i$ must divide $m$.  Because $m$ is a proper divisor of
$a$, this contradicts the minimal property of $a$.  [{\bf11}, vol. 2,
p. 402]

Repeated use of (1) replaces each composite factor in a product by a
product of irreducible factors.  Thus one representation as a product
of irreducibles is possible.  Two such representations, involving
different sets of irreducibles are not possible, as (2) shows.  Note
that one irreducible positive integer divides another only if they are
equal since the only divisors of each are $1$ and itself.

\vskip.5\bls{\noindent\bf3. A variant of Euclid's argument.}  

A set is said to be {\it linearly ordered\/}, or to be a {\it chain},
if, given any pair, $a$, $b$, of distinct members, $a\not=b$, one is
greater than the other: $a<b$ or $b<a$.  A set is said to be {\it well
ordered\/} if it is linearly ordered and if each of its nonvoid
subsets contains a {\it least\/} element: one less than all other
members of the subset.  In addition to the set of positive integers,
every finite linearly ordered set is well ordered.  

No set of rational numbers, with their usual order, that contains the
reciprocals of all the positive integers is well ordered.  Examples:
all the rationals, the real numbers, and complex numbers.

The argument we consider here involves a few simple concepts from
twentieth century algebra: rings, semigroups, partially ordered sets,
and chains of divisors.

If $\bb I$ is the ring of integers, $\bb I'$ the set of nonzero
integers, and $\bb S$ the set of positive integers, the last two are
examples of {\it semigroups} whose products are associative and
commutative, and each contains the identity element $1$ for the
products in $\bb I$.  Since $\bb I$ is an integral domain, these
semigroups are {\it cancellative\/}: $xy=xz$ implies $y=z$.

In any ring or semigroup, we define {\it $a$ divides $b$}, or $a\div
b$, to mean the ring or semigroup contains a member $c$ satisfying
$ac=b$.

The difference between a group and a semigroup is the former contains
the reciprocal of {\it each\/} of its members.  In the two semigroups
here, only $1$, in $\bb S$, and only $\pm1$, in $\bb I'$, divides {\it
every\/} member of the semigroup, including the identity element.
Members of rings or semigroups with this property are called {\it
units}.

Possibly more esoteric than things mentioned so far, although surely
simple as required here, are the MacKenzie {\it co-ideals\/} in a
semigroup [{\bf15}]:

\noindent {\smc DEFINITION. \it If $S$ is a commutative semigroup and
$a\in S$, the subset $\{a\}$ of $S$ containing all divisors of $a$ is
called a co-ideal of $S$.}

The relation $a\div b$ is {\it reflexive\/}, since $a\div a$ for all
$a\in S$, and {\it transitive\/} since $a\div b$ and $b\div c$ imply
$a\div c$, for all $a,b,c\in S$; thus the relation $a\div b$ shares
two of the four defining properties of a `weak' linear order relation, $a\le
b$.  In the {\it special case\/} of the positive integers, it is also
{\it antisymmetric\/}, since for each $a$, $b$,
\empty$$\displaylines{
a\div b \hbox{ and } b\div a \words{imply} a=b, \cr
\hbox{which corresponds to the linear order property}\hfill \cr
a\le b \hbox{ and } b\le a \words{imply} a=b. \cr}
\empty$$
The only other property of a linear order relation that $a\div b$
lacks in general is, given any $a,b \in S$,
\empty$$
\hbox{either $a\le b$ or $b\le a$,\Quad if $a\not=b$.}
\empty$$

A co-ideal $\{a\}$, is a finite subset of $\bb S$, {\it partially
ordered\/} by the divisor relation $x\div y$ in $\bb S$, that we
may also write as $x\le y$ using customary notation for a partial
order relation.  

As partially ordered sets, co-ideals can be represented as the union
of chains by a theorem of Dilworth [{\bf6}, theorem 1.1]; for the
finite case, sufficient here, see Tverberg [{\bf18}].
   
Each $x$ in such a chain is a proper divisor, $x<y$, of each $y$ above
it in the chain.  These chains, being finite linearly ordered sets are
well ordered, and $S$ is said to satisfy the {\it divisor chain condition}.

In terminology for partially ordered sets, $a$ is the universal
(greatest) element of the co-ideal $\{a\}$ generated by $a$, and the
identity element, $1$, for products in the ring containing $\bb S$ is
the null (least) element of $\{a\}$.  In the following example, $a$ is
a square free integer and the co-ideal $\{a\}$ is a boolean algebra.
In less special cases, co-ideals can have more complicated
structures. [Cf.\tsp\St10, below.]

{\noindent\smc EXAMPLE:} The co-ideal of the divisors of 30 has a total
of $2^3$ members and is the union of the following three chains.
\empty$$
1,\;2,\;6,\;30;\Quad 1,\;3,\;15,\;30;\Quad 1,\;5,\;10,\;30.
\empty$$
The {\it atoms\/} of a co-ideal, \ie, members covering only the null
element, are {\it irreducible\/} elements of $\bb S$ since $1$ is their
only proper divisor.

In a partially ordered set if either $a\le b$ or $b\le a$, $a$ and $b$
are called a {\it comparable pair}, otherwise $a$ and $b$ are called
{\it noncomparable}.  According to Dilworth's theorem involved here, a
partially ordered set $P$ can be decomposed into the union of $k$
chains provided each subset of $P$ with $k+1$ elements contains
a comparable pair.

The atoms in $\{a\}$ have the maximal number of noncomparable elements
for $\{a\}$ and the number of chains is the same as the number of
atoms.

Euclid need not have assumed well ordering of the positive integers
themselves was essential for their factorization into products of
irreducibles.  Here we have found well ordering where it is needed, in
chains of divisors of integers, which is just where Euclid used it.

Euclid's criterion (2), or something effectively equivalent to it, is
essential for uniqueness of the representation of positive integers as
products of irreducibles and is what rules out the possibility of
obtaining distinct products that are equal, say,
\empty$$
p_1p_2\cdots p_k = q_1q_2\cdots q_l, \Quad p_i, q_j \hbox{ irreducibles,}
\empty$$
where the intersection of the set of $p$s and the set of $q$s is empty.

\break
\vskip.5\bls{\noindent\bf4. Liouville rings of arithmetic functions.} $(F)$

If $F$ is the {\it field\/} of rational, real, or complex numbers, we
consider rings of arithmetic functions, whose elements are sequences
of members of $F$ indexed by positive integers, and prove unique
factorization for this family of rings.  In \St9 we consider sequences
of integers.

Let $\bb R$ be a commutative and associative ring defined as follows,
see Dickson [{\bf8}], Bell [{\bf2}], and Apostol [{\bf1}]:
  
     \bgroup \parindent=0pt \parskip.1\bls \def\item{\noindent\Quad}
\item{(1)}
The elements of $\bb R$ are all sequences, indexed by the positive
integers, whose members are elements of $F$. We write 
the sequences in the form
$\alpha(n)$ where $n\in\bb S$, $\alpha$ being an arbitrary function
mapping the positive integers into the field $F$.\par
\item{(2)}
Addition in $\bb R$ is defined simply as
$(\alpha+\beta)(n)=\alpha(n)+\beta(n)$.\par
\item{(3)}
The zero of $\bb R$, \ie, the identity element for addition, is the
con\-stant sequence $\omega(n)=0$, for all $n\in\bb S$.\par
\item{(4)}
Products in $\bb R$ are defined by the finite sums
\empty$$
(\alpha\ast\beta)(n) = \sum_{d\delta=n} \alpha(d)\beta(\delta),
\empty$$
where $d$, $\delta$ run through the positive integers satisfying $d\delta=n$.\par
\item{(5)}
The identity element for products is the sequence defined by\par\vskip-1\bls
\empty$$
\vep(n)=1 \hbox{ for } n=1 \words{and} \vep(n)=0 \hbox{ for } n>1.
\empty$$\par        \egroup

{\noindent Each $\bb R$ is an integral domain:  if
$\alpha,\beta\not=\omega$,  then 
$a, b$ exist with
$\alpha(a)\beta(b)\not=0$,
and $\alpha(n)\beta(m)=0$ for $n<a$, $m<b$; hence $(\alpha\ast\beta)(ab)=
\alpha(a)\beta(b)\not=0$, \ie, products of nonzero elements are not zero.}

Familiar members of each of these rings are Euler's $\vph$-function,
the number of positive divisors of $n$, the sum of those divisors,
etc.

The units in these rings are the divisors of $\vep$.
To find all units we suppose $\alpha\in\bb R$,
$(\alpha\ast\alpha')(n)=1$ for $n=1$, $(\alpha\ast\alpha')(n)=0$ for
$n>1$, and solve for $\alpha'$.  Now with $\alpha(1)\in F$, and
\empty$$ 
\alpha(1)\alpha'(1)=1,
\Space\alpha(1)\alpha'(2)+\alpha(2)\alpha'(1)=0,
\Space\alpha(1)\alpha'(3)+\cdots=0,\dots,
\empty$$ 
we {\it must\/} assume $\alpha(1)\not=0$; this infinite system of
equations can then be solved recursively  by induction
for $\alpha'(1)$, $\alpha'(2)$, $\alpha'(3),\dots$.  
The solution is unique since the coefficient matrix of the first $n$
equations is triangular and all elements on the diagonal  are
$\alpha(1)$.  See Bell [{\bf3}].

Thus all functions $\alpha$ with $\alpha(1)\not=0$ are {\it units\/}
in $\bb R$, and {\it only\/} these.  Observe that the sequences
$\alpha$ in $\bb R$ that are units in $\bb R$ have initial terms
$\alpha(1)$ that are units in the underlying field $F$: every nonzero
element $d$ in a field is a unit since it divides $1$, its cofactor
being its reciprocal: $dd^{-1}=1$.  In $\bb R$, all multiplicative
functions have $\alpha(1)=1$. Therefore the most familiar members of
$\bb R$, being units, are of negligible interest when it comes to
unique factorization.

Consider the nonzero nonunits, $\alpha(1)=0$, in $\bb R$; more
precisely, for some $a>1$, $\alpha(n)=0$ for $n<a$ and
$\alpha(a)\not=0$.  If $\alpha$ and $\beta$ are such nonunits,
$(\alpha\ast\beta)(1) = \alpha(1)\beta(1)=0$.  Also
$(\alpha\ast\beta)(p)= \alpha(1)\beta(p) + \alpha(p)\beta(1)= 0$ for
{\it every\/} prime number $p$.  {\smc CONTRAPOSITIVE}: Every nonzero
nonunit that is {\it not\/} zero for one or more prime numbers $p$
cannot be expressed as a product of two nonzero nonunits; it is an
{\it irreducible element\/}.

Among irreducible elements of $\bb R$: the characteristic function
$\chi$ of any set of prime numbers, \ie, $\chi(n)$ is $1$ if $n$
prime, zero otherwise; $(\pi(n))^2$ or, if $F$ is the real field,
$\log n$, $\log(n!)$, $\log(n^n)$, etc.

Let us call $a$ the {\it rank\/} of $\alpha$ if for $a>1$,
$\alpha(a)\not=0$, and $\alpha(n)=0$ when $n<a$.  If $\beta$ has rank
$b$, then the rank of $\alpha\ast\beta$ is the least $n$ for which 
$(a\ast b)(n)\not=0$.  Since $\alpha(a)\beta(b)\not=0$, and
$\alpha(d)\beta(\delta)=0$ if $d<a$ or $\delta<b$, it follows that 
\empty$$ 
(\alpha\ast\beta)(ab) = \sum_{d\delta=ab} \alpha(d)\beta(\delta) =
\alpha(a)\beta(b)
\empty$$
and the rank of the Liouville product $\alpha\ast\beta$ is $ab$.

If we extend this definition of {\it rank\/} so that units have rank
$1$, it follows that the rank of a product $\alpha\ast\beta$, for
$\alpha$, $\beta$ nonzero, is the product of the ranks of its factors.
Except for the zero element $\omega$, every sequence in $\bb R$ has
finite rank.

Suppose $\alpha$ and $\beta$ have the same rank $q$, that is
$\alpha(q)\not=0$, $\beta(q)\not=0$, and $\alpha(r)=\beta(r)=0$ for
$1\le r<q$.  It follows that $\alpha$, $\beta$ differ by a unit
factor:  assume $\alpha=\beta\ast\gamma$.  By hypothesis  
\empty$$ 
\alpha(q) = \sum_{rs=q} \beta(r)\gamma(s) = \beta(1)\gamma(q) +\cdots+
\beta(q)\gamma(1) = \beta(q)\gamma(1).
\empty$$ 
But $\alpha(q)\not=0$ and $\beta(q)\not=0$ imply $\gamma(1)\not=0$.
If $q$ is a prime, $\alpha$ and $\beta$ are irreducible, $q$ being the
first prime for which both are not zero.  

\vskip.5\bls{\noindent\bf5.  Co-ideals in  Liouville rings.} $(F)$

In Liouville rings we define division and divisors as before: {\it
$\alpha(n)$ divides $\beta(n)$\/} means $\alpha\ast\gamma = \beta$ for
some $\gamma(n)\in\bb R$.  In these cases we write $\bb T$ for the
semigroup of nonzero members of $\bb R$.

Suppose $\{\alpha\}$ is the co-ideal of an arbitrary member
$\alpha(n)$ of $\bb T$; when $\alpha,\beta\in\bb T$, if $\alpha|\beta$
then $\{\alpha\}\subseteq\{\beta\}$, and conversely.  If
$\{\alpha\}\subset\{\beta\}$, we call $\alpha$ a {\it proper
divisor\/} of $\beta$.  If $\alpha\div\beta$ {\it{and}}
$\beta\div\alpha$, then $\{\alpha\}=\{\beta\}$ and conversely.  In the
latter case, we say $\alpha$ and $\beta$ are {\it associates} and for
this relation we write $\alpha\sim\beta$.

If $\alpha(n)$ is a unit, $\{\alpha\}$ contains all the units in $\bb
T$, since every unit, by virtue of dividing every member of $\bb T$,
divides $\alpha$.

Clearly if $\delta$ is a unit and if $\alpha\ast\delta=\beta$, then
$\alpha|\beta$ and, since $\delta$ is a unit, we also have
$\beta|\alpha$ because $\beta\ast\delta'=\alpha$ where $\delta'$ is
the reciprocal of $\delta$; \ie, $\delta\ast\delta'=\vep$. Thus if
$\alpha$ and $\beta$ differ by a unit factor, $\alpha\sim\beta$.

The relation $\alpha\sim\beta$ is at the very least an equivalence
relation, \ie, it is reflexive, transitive, and symmetric:

\item{}
Reflexivity:  $\alpha=\vep\ast\alpha$  ($\vep$ is a unit).
\item{}
Transitivity: $\alpha=\delta_1\ast\beta$ and $\beta=\delta_2\ast\gamma$
imply $\alpha=\delta_1\ast\delta_2\ast\gamma$ ($\delta_1\ast\delta_2$ is a
unit, if $\delta_1$ and $\delta_2$ are units).
\item{}
Symmetry: $\alpha\sim\beta$ implies $\beta\sim\alpha$: if $\delta$ is
a unit and $\alpha\ast\delta=\beta$, then $\beta\ast\delta'=\alpha$ if
$\delta\ast\delta'=\vep$.
\item{}
Congruence: $\alpha_1\sim\alpha_2$ and $\beta_1\sim\beta_2$ imply
$\alpha_1\ast\beta_1\sim\alpha_2\ast\beta_2$, since 
$\alpha_1=\delta_1\ast\alpha_2$ and $\beta_1=\delta_2\ast\beta_2$ imply
$\alpha_1\ast\beta_1=(\delta_1\ast\delta_2)\ast\alpha_2\ast\beta_2$.\par

{\noindent}With the last property included, the equivalence relation
shares an additional well-known property with congruence relations.

For $r\ge1$ let $\nu_r(n) = 1$ if $n=r$ and $\nu_r(n)=0$ if $n\not=r$.
Since the functions in each class of associates have the same rank, each
class contains one of the functions $\nu_r$; and if $p$ and $q$ are
equal or distinct positive integers, $\nu_p\ast\nu_q=\nu_{pq}$, \ie,
$(\nu_p\ast\nu_q)(n) = \nu_{pq}(n)$ for $n\ge1$.

\vskip.5\bls{\noindent\bf6.  Reduced semigroups of Liouville rings.} $(F)$

The elements of $\bb T$ are all the nonzero members of $\bb R$,
including, along with each member $\alpha$, all of its associates.
Let $\obar\alpha$ be the equivalence class whose members are $\alpha$
and its associates.  Thus the semigroup $\bb T$ is the union of all
the classes $\obar\alpha$ for $\alpha\in\bb T$, and if $\alpha$ and
$\beta$ are not associates, the intersection $\obar\alpha\cap\obar\beta$ is
empty.

For simplicity, we resort to a {\it homomorphism}, \ie, a many-to-one
mapping from a larger set to a smaller set; we need one that maps
Liouville products in $\bb T$ into Liouville products in the
collection $\obar{\bb T}$ of its equivalence classes.  Thus we
can transfer attention to a semigroup whose {\it members\/} are the
equivalence classes of $\bb T$.  We define Liouville products in the
semigroup $\obar{\bb T}$ of equivalence classes as follows.

Given classes $\obar\alpha$ and $\obar\beta$, we take one element from
each class, say $\alpha_1\in\obar\alpha$ and $\beta_1\in\obar\beta$,
and form their product $\alpha_1\ast\beta_1$.

According to the congruence-condition in the previous section, the
product $\alpha_1\ast\beta_1$ is in the same class with every other
product of the form $\alpha_2\ast\beta_2$ where $\alpha_2$ is {\it
any\/} member of $\obar\alpha$, and $\beta_2$ is {\it any\/} member of
$\obar\beta$.  Since the class containing these products depends only
on the classes used and not on the members multiplied, we define
\empty$$
\obar\alpha\ast\obar\beta \words{to be} \obar{\alpha\ast\beta}.
\empty$$
By shifting attention from  $\bb T$ to $\bb{\obar T}$  we obtain a
major simplification: note, for one example, the equivalence class
containing all the units of $\bb T$ is the identity element for the
Liouville product in $\obar{\bb T}$.

$\obar{\bb T}$ is called the {\it reduced semigroup\/} of $\bb T$.

{\noindent\bf7.  Isomorphism of the reduced semigroups $\bb S$ and
$\bb{\obar T}$.} $(F)$

We have seen the associates of $\alpha$ all have the same rank.  So
there is a one-one correspondence between classes of associates and
the common rank of their members: $\obar\alpha\leftrightarrow a$, 
$\obar\alpha\in\bb{\obar T}$, $a\in\bb S$.

The co-ideals $\{a\}$ and $\{\obar{\alpha}\}$ are {\it isomorphic}
as partially ordered sets, \ie, if $\obar\alpha$ in the latter co-ideal
is paired with its rank $a$ in the former co-ideal, the divisor
relations correspond since $\obar{\beta}\div\obar{\alpha}$ in
$\obar{\bb T}$ if and only if $b\div a$ in $\bb S$.  Therefore the
partially ordered sets $\{a\}$ and $\{\obar\alpha\}$ are identical,
except for the alphabets used to describe them.

{\noindent\bf8. Factorization into irreducibles and its uniqueness.} $(F)$

For $\obar{\bb T}$ we argue as for $\bb S$, with the result that
factorization into finitely many irreducible equivalence classes is
possible and unique.

Each irreducible equivalence class can be represented by an arbitrary
member of the class and alternate representatives may be chosen if in
each instance a compensating unit factor is introduced as well;
moreover, no matter how many alternate choices are made, the resulting
product of units is itself a unit.

Thus each composite function $\alpha$, can be written as the product
of a single unit and finitely many irreducible factors.

Hence factorization of the {\it elements\/} of $\bb R$ takes the form
\empty$$
\delta\,\pi_1\,\pi_2\cdots\pi_k,
\empty$$
where $\delta$ is a unit and $\pi_1, \pi_2,\dots, \pi_k$ are irreducible
members of $\bb R$, a representation that is unique, except for the order
of the factors and replacement by associates of the factors present. 

{\noindent\bf9. Special Case: The subring of sequences of integers.} $(\bb I)$

The first proof for unique factorization in a ring of arithmetic
functions was carried out for the field of complex numbers, a special
case of the family of fields $(F)$ described above in \St4; see
Cashwell and Everett [{\bf4}].  These authors conjectured unique
factorization should also hold for the ring of integer valued
arithmetic functions and proved it soon after [{\bf5}].  In both
cases $F$ and $I$ their methods were transfinite.

We now consider the case of the Liouville ring of arithmetic functions
with integral values and encounter no significant differences from our 
elementary proofs when the functions have values in a field.

For this special case, we observe first that integers form a subring
of each field $F$ considered before, and the ring of sequences of
integers is a subring of each of the corresponding rings $\bb R$.  The
argument presented in sections 4 through 8 can be carried out simply
by restricting attention to sequences of integers.

There are two differences between the cases based on $F$ and $\bb I$.

One is the paucity of units: for sequences of integers units have
$\alpha(1)=\pm1$ instead of $\alpha(1)\not=0$ as for fields.  This
difference is negligible, the homomorphism simply wipes it out.

The difference in what qualifies as a unit entails the other difference:
there are more nonunits, indeed there is a new family of
irreducible elements.  Because the defining condition for units now is
$|\alpha(1)|=1$ instead of $\alpha(1)\not=0$, the two kinds of
nonunits are distinguished since they satisfy {\it either
$\alpha(1)=0$ or $|\alpha(1)| \ge2$}.

The case $\alpha(1)=0$ is common to the earlier cases and the integral
case, and presents no novelty.

The case $|\alpha(1)|\ge2$ is even simpler.  Observe that if
$\alpha(n)=\beta(n)\ast\gamma(n)$ for all $n\ge1$ and if
$|\alpha(1)|=|\beta(1)|$, then $\gamma(1)=\pm1$ since
$|\alpha(1)|=|\beta(1)|\cdot|\gamma(1)|$.  These nonunit sequences
have rank one, in the terminology of the previous cases, and their only
divisors are associates and units; therefore they are irreducibles.

In spite of these differences, the variant of Euclid's argument goes
through as for the fields $F$ considered earlier.

\begingroup  \def\join{\vee}   \def\meet{\wedge}  

{\noindent\bf10. Summary of criteria for unique factorization.}

The Euclidean algorithm ({\sl Elements, \smc VII,\tsp{1},\tsp2}),
often considered in this context, is obtained by repeated subtraction,
which is easy for positive integers because they are (well) ordered
and the division algorithm is available; but not in $\bb R$.

Gauss argues that if a prime number divides neither $a$ nor $b$, it
cannot divide their product; this is a version of Euclid's condition
(2).  Gauss proves it [{\bf13},\tsp{\smc II,\tsp14,\tsp15}] using
congruences, also based on division with remainder and difficult to
achieve without an order relation.

Dirichlet, in his lectures [{\bf9},\tsp\St1.8] avoids {\smc GCD}s and
the division algorithm, and argues as follows.  If $p$ is irreducible,
its only divisors are $1$ and $p$, so $p$ shares no other divisors with
{\it any\/} number.  Dirichlet's observation implies Gauss's version
of Euclid's criterion (2), as an immediate consequence of the
definition of irreducibles.

For the cases considered here the essential hypothesis is the divisor
chain condition, and for this have co-ideals.

For proofs in $\bb I$ and $\bb R$ the following hypotheses are
superfluous: the division algorithm, the Euclidean algorithm, Euclid's
(2), and existence of {\smc GCD}s. The more relevant of these come
from the co-ideals and the definition of irreducibles.

Co-ideals are lattices: We claimed the co-ideal of divisors of 30 is a
boolean algebra because 30 has no square divisors.
              
A boolean algebra is a distributive lattice with unique complements;
its characteristic property is that it is isomorphic to the lattice of
all subsets of some set.  The co-ideal $\{30\}$ is isomorphic to the
collection of all subsets of the set $\{2,3,5\}$, including both the
empty subset and the full set representing the divisors $1$ and $30$.

If in a partially ordered set each pair of elements $x$, $y$ has a
greatest lower bound and a least upper bound in its partial order
relation $x\le y$, the partially ordered set is called a {\it
lattice\/}; the {\smc GLB} and {\smc LUB} are called, respectively,
the {\it meet\/} and {\it join\/} of $x$ and $y$ and are denoted by
$x\wedge y$ and $x\vee y$.

If the partial order relation is the divisor relation $x\div y$,
$x\meet y$ is the {\smc GCD} of $x$ and $y$ and $x\join y$ is their
{\smc LCM}.  See P\'olya and Szeg\H{o} [{\bf16}].

A lattice is {\it complemented\/} if (a) it contains a universal
element and a null element, often denoted by $I$ and $O$,
respectively; and (b) if for each member $x$ there is at least one
member $x'$ with the properties $x\meet x'=O$ and $x\join x'=I$.  (c)
If each $x$ has exactly one complement $x'$, the lattice is said to
have {\it unique complements}.

A lattice is {\it distributive\/} if $ x\vee(y\wedge z) = (x\vee
y)\wedge(x\vee z)$ holds for all its members $x$, $y$, $z$. A
necessary and sufficient condition for distributivity is $x\wedge y =
x\wedge z$ and $x\vee y = x\vee z$ imply $y=z$, Crawley and Dilworth
[{\bf6}, p.\tsp21].  If the partial order is the divisor relation
$x\div y$ both of these assertions are properties of {\smc GCD}s and
{\smc LCM}s.  The necessary and sufficient condition follows
immediately from $(x,y)[x,y] = xy$ and also implies that {\it any\/}
complements that exist are unique.

For the special case of the co-ideal $\{30\}$, the complementary pairs
are $1'=30$, $2'=15$, $3'=10$, $5'=6$, and by symmetry, $30'=1$, etc.
It follows that the co-ideal $\{30\}$ is a boolean algebra.

The co-ideal $\{12\}$ is the union of the following pair of chains
$\langle 1,2,4,12\rangle$ and $\langle 1,3,6,12\rangle$.
Complementary pairs of elements are $1'=12$, and $3'=4$, but not $2$
and $6$ since $2\meet6=2$ and $2\join6=6$. Here $12$ is not square
free; it is divisible by the square of 2.  For more on lattices, see
Crawley and Dilworth [{\bf6},\tsp chapters\tsp1,\tsp2].

For earlier studies of Liouville rings see {\bf4}, {\bf5}, and
{\bf7}; these papers use transfinite methods.  The full force of
Dilworth's theorem does also; it is proved using Zorn's lemma.

For the relation between arithmetic functions and Dirichlet series,
see P{\'o}lya and Szeg{\H o}, [{\bf16}].  Cashwell and Everett's
articles are based on the fact that Dirichlet series may be
represented as power series in infinitely many variables, one variable
for each positive prime number.  In 1913 Harald Bohr published a study
of the analytic properties of this pair of series representations in
{\sl G{\"o}ttinger Nachrichten}, pages 444\dash 488 [Reprinted as Item
A\tsp9, in volume one of Bohr's {\sl Collected Mathematical Works},
K{\o}benhavn (1952)].

The product $\ast$ in $\bb R$ is often called Dirichlet convolution:
see Apostol [{\bf1}, \St 2.14] and Popken [{\bf17}].
              \endgroup

  \def\ref#1:#2\eref{\par\item{{\bf#1\Space}}#2} \overfullrule=0pt

{\vskip .5\bls \parindent=12pt \frenchspacing \raggedright
\centerline{\bf References}\parskip 2.5pt

\ref 1: Tom M. Apostol, {\sl Introduction to Analytic Number Theory},\break
Springer (1976).  Chapter 2.\eref

\ref 2: Eric Temple Bell, Abstract \#4, {\sl Bulletin of the American 
Mathematical Society}, {\bf19} (1913), 166\dash167.\eref

\ref 3: Eric Temple Bell, The reciprocal of a numerical function, {\sl
Tohoku Mathematical Journal}, {\bf 43} (1937), 77\dash78.\eref

\ref 4: E.\tsp D. Cashwell and C.\tsp J. Everett, The ring of
number-theoretic functions, {\sl Pacific Journal of Mathematics},
{\bf9} (1959), 975\dash985.  Math.\tsp Revs.  {\bf21} (1960),
p.\tsp 1334, \#7226  (de Bruijn). \eref

\ref 5: E.\tsp D. Cashwell and C.\tsp J. Everett, Formal power series,
{\sl Pacific Journal of Mathematics}, {\bf13} (1963), 45\dash64.
Math.\tsp Revs. {\bf27} \break (1964), p.\tsp 1101, \#5786. \eref

\ref 6: Peter Crawley and Robert P. Dilworth, {\sl Algebraic Theory of
Lattices}, Prentice-Hall (1973).\eref

\ref 7: Don Deckard and L. K. Durst, Unique factorization in power
series rings and semigroups. {\sl Pacific Journal of Mathematics},
{\bf16} (1966), 239\dash242. Math.\tsp Revs. {\bf32} (1966), p.\tsp
410, \#2439.  \eref

\ref 8: Leonard Eugene Dickson, {\sl History of the Theory of
Numbers}, three volumes (1919), (1920), (1923).  Carnegie Institution
of Washington, publication 256. Reprinted Hafner (1934), Chelsea
(1942).  Vol. 1, pp. 285, 286.\eref

\ref 9: P. G. L. Dirichlet, {\sl Lectures on Number Theory}, History
of Mathematics Sources {\bf16}, American Mathematical Society/London
Mathematical Society (1999).\eref

\ref 10: Harold M. Edwards, {\sl Fermat's Last Theorem}, Springer
(1977). \S\St4.1, 4.2.\eref

\ref 11: Euclid, {\sl The Thirteen Books of The Elements}, translated
with commentary by Sir Thomas L. Heath.  Cambridge University Press,
second edition (1926). Reprinted Dover (1956) and later.\eref

\ref 12: G. H. Hardy and E. M. Wright, {\sl An Introduction to the 
Theory of Numbers}, Oxford (1938).  Fifth edition (1979). Chapter
14.\eref

\ref 13: Carl Friedrich Gauss, {\sl Disquisitiones Arithmetic\ae},
English,\break Springer (1986).\eref

\ref 14: Nathan Jacobson, {\sl Basic Algebra I}, W. H. Freeman, 2nd.
edition (1985). \St2.14. \eref

\ref 15: Robert\tsp E.\tsp MacKenzie, Commutative semigroups.  {\sl
Duke Mathematical Journal}, {\bf21} (1959), 975\dash 985.\eref

\ref 16: George P\'olya and Gabor Szeg\H{o}, {\sl Problems and
Theorems in\break Analysis}, Springer, two volumes. First German ed. (1925),
English (1978), (1976), reprinted (1998). Part VIII, \S\S\tsp1.4,\tsp1.5.\eref

\ref 17: Jan Popken, On convolutions in number theory, {\sl
Koninklijke Nederlandes Akademie van Wetenschappen\/}: Proceedings
Series A. Mathematical Sciences (Amsterdam), {\bf58} = {\sl
Indagationes Mathematic\ae}, {\bf17} (1955), 10\dash15.  \eref

\ref 18: Helge Tverberg, On Dilworth's decomposition theorem for partially
ordered sets. {\sl Journal of Combinatorial Theory}, {\bf3} (1967),
305\dash306.\eref

\hfil

\item{\copyright}May 2001, Lincoln Durst. \Quad \tt g\_lkd@ams.org}

\bye